\documentclass[reqno, 11pt]{amsart}
\usepackage[colorlinks]{hyperref}
\usepackage{url}
\newtheorem{theorem}{Theorem}[section]
\newtheorem{lemma}[theorem]{Lemma}
\newtheorem{proposition}[theorem]{Proposition}
\newtheorem{corollary}[theorem]{Corollary}

\theoremstyle{definition}
\newtheorem{definition}[theorem]{Definition}

\theoremstyle{remark}
\newtheorem{remark}[theorem]{Remark}

\numberwithin{equation}{section}

\usepackage{bbm}
\usepackage{euscript}
\usepackage{pb-diagram}
\usepackage{lamsarrow}
\usepackage{pb-lams}
\usepackage{amsmath}
\usepackage{amsthm}
\usepackage{epsfig}
\usepackage{amssymb}
\usepackage{pifont}
\usepackage{amsbsy}

\newskip\aline \newskip\halfaline
\def\skipaline{\vskip\aline}

\def\qedbox{$\rlap{$\sqcap$}\sqcup$}
\def\qed{\nobreak\hfill\penalty250 \hbox{}\nobreak\hfill\qedbox\skipaline}

\def\proofend{\eqno{\mbox{\qedbox}}}

\newcommand\bR{{\mathbb R}}
\newcommand\bS{{\mathbb S}}

\newcommand\bZ{{\mathbb Z}}

\newcommand{\bn}{\boldsymbol{n}}
\newcommand{\bp}{\boldsymbol{p}}

\newcommand{\bu}{\boldsymbol{u}}

\newcommand{\bsA}{\boldsymbol{A}}

\newcommand{\bsN}{\boldsymbol{N}}
\newcommand{\bsR}{\boldsymbol{R}}

\newcommand{\btau}{\boldsymbol{\tau}}

\newcommand{\si}{{\sigma}}

\newcommand{\ve}{{\varepsilon}}

\newcommand{\vfi}{{\varphi}}

\newcommand{\eC}{\EuScript{C}}

\newcommand{\eL}{\EuScript{L}}
\newcommand{\eM}{\EuScript{M}}
\newcommand{\eN}{\EuScript{N}}
\newcommand{\eO}{\EuScript{O}}

\newcommand{\eS}{\EuScript{S}}

\newcommand{\eU}{\EuScript{U}}
\newcommand{\eV}{\EuScript{V}}

\newcommand{\eX}{\EuScript{X}}

\newcommand{\ra}{\rightarrow}
\newcommand{\hra}{\hookrightarrow}

\newcommand{\lan}{\langle}
\newcommand{\ran}{\rangle}

\def\inpr{\mathbin{\hbox to 6pt{\vrule height0.4pt width5pt depth0pt \kern-.4pt \vrule height6pt width0.4pt depth0pt\hss}}}

\newcommand{\pa}{\partial}

\DeclareMathOperator{\cl}{\boldsymbol{cl}}

\DeclareMathOperator{\Cr}{\mathbf{Cr}}

\begin{document}

\title{On the total curvature  of semialgebraic graphs}
\date{Version 1. Started: June 6, 2008. Completed: June 18,2008.}

\author{Liviu I. Nicolaescu}

\address{Department of Mathematics, University of Notre Dame, Notre Dame, IN 46556-4618.}
\email{nicolaescu.1@nd.edu}
\urladdr{http://www.nd.edu/~lnicolae/}

\begin{abstract}   We  define the total curvature  of a semialgebraic graph  $\Gamma\subset \bR^3$   as an integral $K(\Gamma)=\int_\Gamma d\mu$,  where $\mu$ is  a  certain Borel measure    completely determined by the local extrinsic geometry of $\Gamma$.  We prove that it  satisfies the Chern-Lashof inequality  $K(\Gamma)\geq b(\Gamma)$, where $b(\Gamma)=b_0(\Gamma)+b_1(\Gamma)$, and we  completely characterize    those graphs for which we have equality.   We also prove the following unknottedness result:  if $\Gamma\subset \bR^3$ is homeomorphic to the suspension of an $n$-point set, and satisfies  the inequality   $K(\Gamma) <2+b(\Gamma)$,  then  $\Gamma$ is unknotted.  Moreover, we  describe a simple  planar graph  $G$  such that for any $\varepsilon>0$ there exists a knotted  semialgebraic embedding  $\Gamma$ of $G$ in $\bR^3$ satisfying  $K(\Gamma)<\varepsilon +b(\Gamma)$.
\end{abstract}

\maketitle

\tableofcontents

\section*{Introduction}
\addtocontents{section}{Introduction}
The total curvature of  a simple closed $C^2$-curve $\Gamma$ in $\bR^3$  is the quantity
\[
K(\Gamma)=\frac{1}{\pi}\int_\Gamma|k(s)\,|ds|,
\]
where $k(s)$ denotes the curvature function of $\Gamma$ and $|ds|$ denotes the arc-length along $C$. In 1929 W. Fenchel \cite{Fen} proved that for  any  such curve $\Gamma$ we have the inequality
\begin{equation*}
K(\Gamma)\geq 2,
\tag{$\mathrm{F}$}
\label{tag: f}
\end{equation*}
with equality  if and only if $C$ is  a planar convex curve.

Two decades later, I. F\'{a}ry \cite{Fa} and J. Milnor \cite{Mil} gave   probabilistic interpretations of the total curvature.  Milnor's interpretation goes as follows.

Any unit vector $\bu\in\bR^3$ defines a linear function $h_{\bu}:\bR^3\ra \bR$, $x\mapsto (\bu,x)$, where $(-,-)$ denotes the inner product in $\bR^3$. For a generic $\bu$, the restriction of $h_{\bu}$ to $\Gamma$ is a Morse function. We denote by $w_\Gamma(\bu)$ the number of critical points of  this function.  Then
\[
K(\Gamma)=\frac{1}{{\rm area}\,(S^2)} \int_{S^2}w_\Gamma(\bu)\,|d\bu|,
\]
where $S^2$ denotes the unit sphere in $\bR^3$ and $|du|$ the Euclidean area density on $S^2$.  Since any function on $\Gamma$ has at least two critical points (a minimum and a maximum) the  inequality (\ref{tag: f}) is obvious. Moreover, they show  that if $K(\Gamma)$ is not too large, then  $\Gamma$ cannot be knotted. More precisely, if $K(\Gamma) <4$ then $\Gamma$ cannot be knotted.

Soon after, in 1957,  Chern and Lashof \cite{CL}    proved  higher dimensional   generalizations of the  results of  Fenchel, F\'{a}ry and Milnor.   Fix a compact $k$-dimensional  submanifold  $\Gamma\subset \bR^{n+1}$.  Again, any unit vector $\bu\in \bR^{n+1}$ defines a linear function $h_{\bu}$ on $\bR^{n+1}$. For generic $\bu$ restriction of $h_{\bu}$ to $\Gamma$ is a Morse function. We denote by $w_\Gamma(\bu)$ the number of its critical points.  Observe that if $M_{\bu}(t)$ denotes the Morse polynomial of $h_{\bu}|_\Gamma$ then $w_\Gamma(t)=M_{\bu}(t)|_{t=1}$. We set
\[
K(\Gamma)=\frac{1}{{\rm area}\,(S^n)}\int_{S^n} w_\Gamma(\bu)\,|du|,
\]
where $S^n$ denotes the unit sphere in $\bR^{n+1}$.

The  Morse inequalities  imply that   $w_\Gamma(\bu)\geq \sum_{j=0}^kb_j(\Gamma)$  for generic $\bu$, where $b_j(\Gamma)$ are the Betti numbers of $\Gamma$.  In particular, we obtain the Chern-Lashof inequality
\begin{equation*}
K(\Gamma)\geq \sum_{j=0}^k b_j(\Gamma).
\tag{$\mathrm{CL}$}
\label{tag: cl}
\end{equation*}
Chern and Lashof proved that, much like in the case of curves, the     quantity $K(\Gamma)$ can be expressed as an integral
\[
K(\Gamma) = \int_\Gamma  \rho_\Gamma (x) |dA_\Gamma(x)|,
\]
where $\rho_\Gamma(x)$ can be explicitly computed from the second fundamental form of the embedding $\Gamma\hra \bR^{n+1}$, and $|dA_\Gamma|$ is the Euclidean area density on $\Gamma$. Additionally, they proved that $K(\Gamma)=2$ if and only if $\Gamma^k$ is a convex hypersurface  of an affine $(k+1)$-dimensional plane in $\bR^{n+1}$. The  embedding $\Gamma^k\hra \bR^{n+1}$ is called \emph{tight} if we have equality in (\ref{tag: cl}).   The subject of tight embeddings continues  to be  an active area of research (see e.g. \cite{Ba,Ku,Kui}).

In this  paper we extend  the Chern-Lashof approach to singular
one dimensional compact semialgebraic subsets of $\Gamma\subset
\bR^3$.   They can be visualized as graphs embedded in some
``tame''\footnote{For example, tameness  would prohibit  ``very
wavy'' edges.} fashion in $\bR^3$. There are several competing
proposals of what should constitute the total curvature of such a
graph (see e.g. \cite{GY, Tan})) but they don't seem to fit the
elegant mold    created by Chern and Lashof. Our approach
addresses precisely this issue and its   uses an approach based on
stratified Morse theory  pioneered  by T. Banchoff \cite{Ba0}  and
N. Kuiper \cite{Kui1} for special cases of  stratified  spaces,
more precisely, PL spaces. Here are the main ideas and results.

Consider  a compact, connected  one-dimensional semialgebraic subset $\Gamma\subset\bR^3$.   We  fix a Whitney stratification of $\Gamma$, i.e., we   fix a finite  subset $V\subset \Gamma$ such that the complement is a  finite disjoint union of $C^2$ arcs.  Then, for a generic $\bu\in S^2$ the restriction of $h_{\bu}$ to $\Gamma$ is a stratified Morse function in the sense of Lazzeri \cite{Laz} and  Goreski-MacPherson \cite{GM}. We denote  by $M_{\bu}(t)$ its stratified Morse polynomial and we set
\[
w_\Gamma(\bu)=M_{\bu}(t)|_{t=1}.
\]
The  stratified Morse inequalities imply that $w_\Gamma(\bu)\geq b_0(\Gamma)+b_1(\Gamma)=1+b_0(\Gamma)$ and we define the total curvature of $\Gamma$ to be
\[
K(\Gamma)=\frac{1}{{\rm area}\,(S^2)}\int_{S^2} w_\Gamma(\bu)\,|du|.
\]
Clearly the total curvature  satisfies the Chern-Lashof inequality (\ref{tag: cl}), and we say that $\Gamma$ is \emph{tight}   if we have equality.

In  Theorem \ref{th: CL} we  give an explicit description of $K(\Gamma)$ in terms of infinitesimal and local invariants of $\Gamma$ which  shows that the total curvature is independent  of the choice of the Whitney stratification. The can be given a characterization  similar in spirit to the approach in  \cite{Mil}. More precisely (see  Corollary \ref{cor: min}) the number $\mu(\Gamma)=\frac{1}{2}\bigl(\,K(\Gamma)+\chi(\Gamma)\,\bigr)$ is equal to the  average number of  local minima of the  family of functions $h_{\bu}|_{\Gamma}$, $\bu\in S^2$.  Following  the terminology in \cite{Mil} we will refer to $\mu(\Gamma)$ as the \emph{crookedness} of $\Gamma$. The Chern-Lashof  inequality can be rephrased as $\mu(\Gamma)\geq 1$.

In  Corollary \ref{cor: GY} we proved that  if the vertices  of  $\Gamma$ have degrees $\leq 3$  then our  integral curvature coincides (up to a multiplicative factor) with the integral curvature recently introduced by Gulliver and Yamada \cite{GY}. In general there does not seem to be a simple relationship between these two notions of integral curvature.

We  also  investigate the structure of one-dimensional tight
semialgebraic   sets.  We observe that $\Gamma$ is tight if and
only if it satisfies  Banchoff's  \emph{two-piece property}: the
intersection of $\Gamma$ with any closed half-space  is either
empty, or connected.  Using this observation we were able to give
a complete description of the tight one dimensional semialgebraic
subsets of $\bR^3$. More precisely, in Theorem \ref{th: tight} we
prove that they are of two types.

\medskip

\noindent $\bullet$ \emph{Type S: Straight.}  In this case all the edges are \emph{straight} line segments. Moreover, there exists   a convex polyhedron (canonically determined by $\Gamma$ such that the following hold (see \cite[Lemma 2.4]{Ku})

\begin{enumerate}

\item[(a)] The  $1$-skeleton of $P$ is contained in $\Gamma$.

\item[(b)] Any vertex $v$ of $\Gamma$   which is not a vertex of $P$  has the property  lies in the convex hull of its neighbors.

\end{enumerate}

\noindent $\bullet$ \emph{Type C: Curved.}  In this case some of the edges of $\Gamma$  have nontrivial curvature. Then $\Gamma$ is contained in  a plane $P$ and there exists a closed convex semialgebraic curve $B\subset P$ with the following properties.

\begin{enumerate}

\item[(a)] $B\subset \Gamma$.

\item[(b)] $\Gamma \setminus B$ is a union of  line segments contained in the region $R$ bounded by $B$.

\item[(c)] The complement of $\Gamma$ in the region bounded by $B$ is a finite union of convex open subsets of the plane $P$.

\end{enumerate}

In particular, this gives a  positive answer to a question   raised at the end of \cite[Sec. 4]{GY}.

We also discuss knottedness issues.    In Theorem \ref{th: unknot} we prove  that if
$\Gamma\subset \bR^3$ is a  semialgebraic subset of $\bR^3$
homeomorphic to the suspension of an $n$-point set, and $\mu(\Gamma)<2$, then $\Gamma$ is isotopic to a planar embedding of this suspension.  The case $n=2$  was first proved by  F\'{a}ry \cite{Fa} and  Milnor \cite{Mil}, while the case $n=3$ was investigated  Gulliver-Yamada \cite{GY} who proved the  unknottedness   under the  more stringent requirement $\mu(\Gamma)<\frac{3}{2}$.

 The situation is dramatically different for
slightly more complicated graph. Consider a graph  which  is homeomorphic to the union of a round circle and two parallel chords.   We show that for every $\ve>0$ there exists a knotted $PL$-embedding $\Gamma_\ve\hra \bR^3$ of this graph such that
\[
\mu(\Gamma_\ve)<1+\ve.
\]
For the reader's convenience, we have included a brief appendix containing    some basic facts about semialgebraic sets  used throughout the paper.

\medskip

\noindent {\bf Notations} In this paper, we will denote by $(-,-)$ the inner product in $\bR^3$, by $|\bullet|$ the corresponding Euclidean norm. For any finite set $S$ we will denote its cardinality by $\# S$.

\section{One dimensional stratified Morse theory}
\setcounter{equation}{0}
Suppose $\Gamma$ is a compact connected $1$-dimensional  semi-algebraic subset  of $\bR^3$.  It can be identified \emph{non canonically} with a graph as follows. We fix a finite subset $V\subset \Gamma$  called the \emph{vertex set} such that the complement $\Gamma\setminus V$ is a finite disjoint  union   of  real analytic,   bounded semialgebraic arcs without self intersections connecting    different points in $V$.  We will refer to these arcs as \emph{open edges} and we will denote  by $E=E(\Gamma)$ the set of open  edges.  Note that   this  definition excludes the existence  edges with  identical endpoints although we allow for multiple edges between two given points.

  A  \emph{semialgebraic graph} is a compact, connected $1$-dimensional   semialgebraic  set together with a choice of vertex set satisfying  the above properties.\footnote{The graph structure is a special Whitney stratification of $\Gamma$.} Clearly, on the same   semialgebraic set we can define multiple structures of semialgebraic sets.

If $\Gamma$ is a semialgebraic  graph with vertex set $V$, then  the degree of a vertex $v$, denoted by $\deg(v)$ is the number of  edges incident to $v$.  Then
\begin{equation}
\chi(\Gamma)=\# V-\# E=\frac{1}{2}\sum_{v\in V}\bigl(\,2-\deg(v)\,\bigr).
\label{eq: eul}
\end{equation}
Since $\Gamma$ is connected  we deduce
\begin{equation}
b_1(\Gamma)= 1-\chi(G)= 1-\# V+\# E.
\label{eq: b1}
\end{equation}
Note that the degree   of a vertex  $p$ can be also defined as the  cardinality of the intersection of $\Gamma$   with a  sphere of sufficiently small radius   centered at $p$. This definition makes sense  even  for  points $p\in \Gamma\setminus V$, and for such points we have $\deg(p)=2$. The equality (\ref{eq: eul}) can be rewritten as
\begin{equation}
\chi(\Gamma)=\sum_{p\in \Gamma} \bigl(\,2-\deg(p)\,\bigr)
\label{eq: eul2}
\end{equation}
For every $p\in \Gamma$ we define $\bsN_p\subset S^2$ as follows.  For $q \in
\Gamma\setminus \{p\}$ denote by $\overrightarrow{\rho_p(u)}$ the  unit vector
$\overrightarrow{\rho_p(q)}:=\frac{1}{|\overrightarrow{pq}|}\overrightarrow{pq}$.
We obtain in this fashion o semialgebraic map
\[
\rho_p:\Gamma\setminus v\ra S^2.
\]
Now set
\[
\btau\in \bsN_p\Longleftrightarrow  \mbox{$\exists$ sequence $(q_k)_{k\geq 1}\subset  \Gamma\setminus \{p\}$},\;\;\btau=\lim_{k\ra \infty}\overrightarrow{\rho_p(q_k)}.
\]
Since $\Gamma$ is semialgebraic, the set $\bsN_p$ is finite for any $p\in \Gamma$. We will refer to the vectors in $\bsN_p$ as the \emph{interior unit tangent vectors} to $\Gamma$ at $v$. The union of  half-lines at $p$  in the directions given by $\btau\in \bsN_p$ is called the \emph{tangent cone} to $\Gamma$ at $p$.

Any unit vector $\bu\in S^2$ defines  a linear map (height function)
\[
h_{\bu}:\bR^3\ra \bR,\;\;h_u(x):=(\bu, x).
\]
 A unit vector $\bu$ is called  $\Gamma$-\emph{nondegenerate}  if the restriction of $h_{\bu}$ to $\Gamma$ is a stratified  Morse function with respect to the  vertex-edge stratification. More precisely,  (see \cite{GM})  this  means that the restriction of $h_{\bu}$ to the interior of any edge has only nondegenerate  critical points and moreover
\[
(\bu,\btau)\neq 0,\;\;\forall v\in V,\;\;\btau\in \bsN_v.
\]
The  complement  in $S^2$ of the  set of $\Gamma$-nondegenerate
vectors is called the \emph{discriminant set} of the semialgebraic
graph and it is denoted by $\Delta_\Gamma$. The discriminant set
clearly  depends on the choice of   graph structure, but one can
prove (see \cite{GM, KaSch} and Lemma \ref{lemma: bertini}) that
$\Delta_\Gamma$ is a  closed semialgebraic subset  of $S^2$ of
dimension $\leq 1$. In particular, most unit vectors $\bu$ are
$\Gamma$-nondegenerate.

\begin{figure}[h]
\centerline{\epsfig{figure=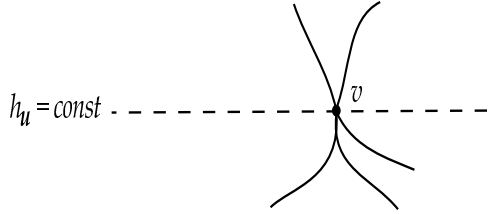,  height=1.1in, width=2.5in}}
\caption{\sl The neighborhood of a critical point.}
\label{fig: 1}
\end{figure}

If $\bu$ is a nondegenerate vector then for every $p\in \Gamma$ and every $\ve>0$ we
set
\begin{equation}
L^\pm_\ve(p,\bu):= \Bigl\{ q\in \Gamma;\;\; |p-q|=\ve,\;\; \pm\bigl(h_{\bu}(q)-h_{\bu}(p)\,\bigr)>0\,\Bigr\},
\label{eq: link}
\end{equation}
\[
d_p^\pm(\bu):=\lim_{\ve\searrow 0} \# L_{\ve}^\pm(p,\bu).
\]
Observe that $d^+_p(\bu)+d^-_p(\bu)=\deg(p)$. In Figure \ref{fig: 1} we have $d^+_v(\bu)=2$, $d_v^-(\bu)=3$, while $\bsN_v$ consists of  four  vectors because two of the edges  are tangent at $v$.

Suppose that $\bu$ is a $\Gamma$-nondegenerate  unit vector. A
\emph{stratified critical point} of $h_{\bu}$ on $\Gamma$ is a
point $p\in \Gamma$  which is either a vertex, or a  critical
point of the restriction of  $h_{\bu}$ to one of the edges. We
denote by $\Cr(\bu)$ the set\footnote{The same unit vector  may be
nondegenerate for several graph structures and the critical sets
corresponding to these graph structures could be different.} of
stratified critical points of $h_{\bu}$. The points in
$\Gamma\setminus \Cr(\bu)$ are  called  \emph{regular points} of
$h_{\bu}$. Observe that if $p$ is a regular point, then
$d^-_p(\bu)=1$.

Note that $\Cr(\bu)$ is a finite subset of $\Gamma$ containing the vertices.  We set $\Cr^1_{\bu}:=\Cr_{\bu}\setminus V$. In other words,  $Cr_{\bu}^1$ consists of the points in the interiors of edges where the tangent vector to the edge is perpendicular to $\bu$.

The set $\Cr^1(\bu)$ decomposes   into a set of local minima $\Cr^1_{\min}(\bu)$, and a set of local maxima $\Cr^1_{\max}(\bu)$. The set of vertices  $V$ further decomposes as $V=V_{\min}(\bu)\cup V^*(\bu)$, where $V_{\min}(u)$ consists of the vertices of $\Gamma$ which are  local minima of $h_{\bu}$, and $V^*(\bu)$ is its complement.  The set $V_{\min}(-\bu)$ is the set of local maxima of $h_{\bu}$, and for this reason we will denote it by $V_{\max}(\bu)$.

Note that  for  every nondegenerate  unit vector $\bu$,  and every $c\in \bR$ the sublevel set $\{h_{\bu}\leq c\}$ is a also a semialgebraic graph, and  the level set  $\{h_{\bu}=c\}$ consists of finitely many points. If a level set $\{h_{\bu}=c\}$ contains only regular points  then  $\{ h_{\bu}\leq c-\ve\} $ is homeomorphic to $\{h_{\bu}\leq c+\ve\}$ for all sufficiently small $\ve$.

If the level set  $\{h_{\bu}=c\}$ contains  the critical points  $p_1,\dotsc, p_k$, then   $\{h_{\bu}\leq c+\ve\}$ is homotopic to the set  obtained from    $\{h_{\bu}\leq c-\ve\}$ by   separately   conning off    each of the sets $L^-_\delta(\bu, p_i)$, $i=1,\dotsc, k$. The points $p_i$ will be the vertices of the added cones.   Here we define the cone over the empty set to be a single point.  For example, if  the level set $\{h_{\bu}=c\}\cap \Gamma$ contains $k$ local minima,   and the sublevel set $\{h_{\bu}\leq c-\ve\}\cap \Gamma$ is connected   for all $\ve>0$ sufficiently small, then  the sublevel set $\{h_{\bu}\leq c+\ve\}\cap \Gamma$ will have $k+1$ connected components for all $\ve>0$ sufficiently small.

To every critical point $p\in \Gamma$ we associate its \emph{Morse polynomial}   $M_{\bu}(t,p)\in \bZ[t]$  according to the rule
\[
M_{\bu}(t,p):= \begin{cases}
1 & \mbox{ if $p$ is a local minimum of $h_{\bu}$} \\
(d_v^-(\bu)-1)t & \mbox{otherwise}.
\end{cases}
\]
Observe that $M_{\bu}(t,p)=0$ if $p$ is a regular point. Let us observe that $M_{\bu}(t, p)$  is the Poincar\'{e} polynomial of the topological pair $(\,{\rm cone}\,(L^-_\ve(p,\bu)), L^-_\ve(p,\bu))$ where ${\rm cone}\,(L^-_\ve(p,\bu))$ denotes the cone over $L^-_\ve(p,\bu)$ and $\ve$ is sufficiently small

The \emph{homological weight} of  the critical point $p\in \Cr(\bu)$ is then defined as the integer
\[
w(p,\bu):=M_{\bu}(t,p)|_{t=1}.
\]
The \emph{Morse polynomial} of $\bu$ is defined as
\[
M_{\bu}(t)=M_{\bu, \Gamma}(t):=\sum_{p\in \Cr(\bu)}M_{\bu}(t,p)=\sum_{p\in \Gamma} M_{\bu}(t,p).
\]
The \emph{homological weight} of $\bu$ is then the integer
\[
w(\bu)=w_\Gamma(\bu):=M_{\bu}(1)=\sum_{p\in \Cr(\bu)}w(p,\bu).
\]
\begin{figure}[h]
\centerline{\epsfig{figure=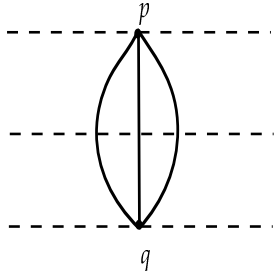,  height=1.4in, width=1.4in}}
\caption{\sl  A stratified Morse function on a $\theta$-graph.}
\label{fig: 4}
\end{figure}

\begin{remark} The homological weight $w_\Gamma(\bu)$  is  in general different from the number of critical points. In  Figure \ref{fig: 4} we have depicted  a stratified Morse function $h_{\bu}$  on  a planar graph. The vector $\bu$ lies in the plane of the graph, it is perpendicular to the dotted lines  and points upwards. This function has precisely two stratified  critical points $p$ and $q$ but its homological weight is $w(\bu)=w(\bu,p)+w(\bu,q)=2+1=3$. Its Morse polynomial is $M_{\bu}(t)=1+2t$. Note that $M_{\bu}(-1)=-1=\chi(\Gamma)$.\qed
\end{remark}
We define a partial order  $\preccurlyeq$ on the  vector space $\bR[t]$ by declaring $P\preccurlyeq Q$ if and only if there exists a polynomial $R$ with \emph{nonnegative coefficients} such that
\[
Q-P=(1+t)R.
\]
The following result  follows immediately  from  the main theorems of stratified Morse theory \cite{GM}.

\begin{theorem}[Morse inequalities]  Denote by $P_\Gamma(t)$ the  Poincar\'e polynomial of $\Gamma$, $P_\Gamma(t)=1+ b_1(\Gamma)t$. Then for every nondegenerate vector  $\bu\in S^2$ we have
\[
M_{\bu}(t)\succcurlyeq P_\Gamma(t).
\]
In particular,
\begin{equation}
M_{\bu}(-1)=P_\Gamma(-1)=\chi(\Gamma)
\label{eq: morse-eul}
\end{equation}
and
\begin{equation}
w_\Gamma(\bu)\geq  1+b_1(\Gamma).
\label{eq: morsein}
\end{equation}\qed
\end{theorem}

\begin{definition} A $\Gamma$-nondegenerate vector $\bu$ is called ($\Gamma$-)\emph{perfect} if $M_{\bu}(t,\Gamma)=P_{\Gamma}(t)$.\qed
\end{definition}

The  situation in Figure \ref{fig: 4}   corresponds  to   a perfect nondegenerate unit vector.

\begin{lemma} A nondegenerate vector $\bu$ is  perfect if and only if  $w_\Gamma(\bu)=P_\Gamma(1)=1+ b_1(\Gamma)$.
\end{lemma}

\proof  Observe that  since $M_{\bu}(t,\Gamma)$ and $P_\Gamma(t)$ are polynomials of  degree $1$, they are completely determined by their values at  two different $t$'s. Since $M_{\bu}(t,\Gamma)|_{t=-1}=P_\Gamma(-1)$ we deduce that  a vector $\bu$ is perfect if and only if $M_{\bu}(1)=P_\Gamma(1)$.\qed

\begin{proposition}  Suppose $\bu$ is  a $\Gamma$-\emph{nondegenerate} vector $\bu\in S^2$. Then the following statements are equivalent.

\smallskip

\noindent (a) The vector $\bu$ is perfect.

\noindent (b) For any $c\in \bR$ the sublevel set $\{h_{\bu}\leq c\}\cap \Gamma$ is connected.
\label{prop: min}
\end{proposition}

\proof   We set $m(\bu)=\min_{p\in \Gamma}h_{\bu}(p)$.  Observe that $M_{\bu}(0)$ is  the number of local  minima of $h_{\bu}$.

\smallskip

\noindent (b) $\Longrightarrow$ (a)   The  sublevel set $\{h_{\bu}\leq m(\bu)\}$ is connected. It coincides with  the set of absolute minima of $h_{\bu}$.  Since $\bu$ is nondegenerate this is a discrete set consisting of a single point $p_0$.       All the sublevel sets $\{h_{\bu}\leq c\}\cap \Gamma$ are connected  so that $h_{\bu}$ cannot have a local minimum other than the absolute minimum.

Indeed, as $c$ increases away from $m(\bu)$, the first time $c$
encounters a critical   value $c_0$  such that the level set
$\{h_{\bu}=c_0\}$ contains local minima, then the number  of
components of $\{h_{\bu}\leq c_0\}\cap \Gamma$ increases by
exactly the number of the local minima.

 This proves  that $1=M_{\bu}(0)=P_\Gamma(0)$. On the other hand, $M_{\bu}(-1)=P_\Gamma(-1)=\chi(\Gamma)$. Since both $M_{\bu}(t)$ and $P_\Gamma(t)$ are polynomials of  degree $1$ we conclude that $M_{\bu}(t)=P_\Gamma(t)=1$.

\smallskip

\noindent (a) $\Longrightarrow$ (b)   Since $\bu$ is perfect we
deduce that $M_{\bu}(0)=P_\Gamma(0)=1$ so that   $h_{\bu}$ has   a
unique local minimum. In particular, the sublevel set
$\{h_{\bu}\leq m(\bu)\}\cap \Gamma$ consists of single point and
thus it is connected. To conclude   run  in reverse the
topological argument in the proof of the implication (b)
$\Longrightarrow$ (a). \qed

\begin{definition} A connected semialgebraic \emph{graph}  $\Gamma\subset \bR^3$ is called \emph{tight} if \emph{all} the  $\Gamma$-nondegenerate vectors $\bu$ are $\Gamma$-perfect, i.e.,
\[
w_\Gamma(\bu)= 1+b_1(\Gamma),\;\;\forall \bu\in S^2\setminus \Delta_\Gamma.
\]
A  compact connected $1$-dimensional semialgebraic \emph{set} $\Gamma\subset \bR^3$ is called \emph{tight} if it is tight for some choice of graph structure.  \qed
\end{definition}

\begin{proposition}  Suppose $\Gamma \subset \bR^3$ is a semialgebraic graph.  Then the following statements are equivalent.

\noindent (a) The graph $\Gamma$ is tight.

\noindent (b) The intersection of $\Gamma$ with {\bf almost  any}
closed half-space is either empty or connected.

\noindent (c) The intersection of $\Gamma$ with {\bf any} closed
half-space is either empty or connected.

\label{prop: tpp}
\end{proposition}

\proof   Observe  first that  the intersection of $\Gamma$  with a closed  half-space is a set of the form $\{h_{\bu}\leq c\}\cap \Gamma$ for some $\bu\in S^2$ and $c\in \bR$.

The implications (c) $\Longrightarrow$ (a), (b) and  (a)
$\Longrightarrow$ (b) follow immediately from Proposition
\ref{prop: min}.  It suffice to prove  only the implication  (b)
$\Longrightarrow$ (c).  We use an argument inspired  by the proof
of \cite[Thm. 3.11]{Kui}.

Denote by $\eU$ the set of  vectors $\bu\in S^2$ such that  for any $c\in \bR$ the sublevel set $\{h_{\bu}\leq c\}$ is  connected.   The  set $\eU$ is dense in $S^2$.

Fix a vector $\bu \in S^2$ and a real number $c$. Then there exists  a sequence  of nondegenerate perfect vectors $\bu_n\in \eU$ and a sequence   of positive real numbers $(r_n)_{n\geq 1}$  satisfying the following properties.

\begin{itemize}

\item $\lim_{n\ra \infty}\bu_n=\bu$.

\item $\lim_{n\ra \infty} r_n=0$

\item $\{h_{\bu}\leq c\}\cap \Gamma\subset \{ h_{\bu_{n+1}}\leq c+r_{n+1}\}\cap \Gamma\subset \{h_{\bu_{n}}\leq c+ r_{n}\}\cap \Gamma$, $\forall n\geq 1$.

\end{itemize}
 The sets $X=\{\, h_{\bu}\leq c\,\}\cap \Gamma$ and $X_n= \{h_{\bu_{n}}\leq c+ r_{n}\}\cap \Gamma$ are closed semialgebraic   subsets of $\Gamma$ satisfying the conditions
 \[
 X\subset X_{n+1}\subset X_n,\;\;\forall n\geq 1\;\;\mbox{and}\;\;  X=\bigcap_{n\geq 1} X_n
 \]
 We conclude that (see \cite[\S 6.6-\S6.8]{Spa})
\[
\check{H}^0(X,\bZ)=\varinjlim_n \check{H}^0(X_n,\bZ),
\]
where $\varinjlim$ denotes the inductive limit and $\check{H}^\bullet(-,\bZ)$ denotes the \v{C}ech homology with integral coefficients. For  semialgebraic sets the \v{C}ech cohomology coincides with the usual singular cohomology, and  all the sets $X_n$ are connected so that $\check{H}^0(X_n,\bZ)=\bZ$ for all $n$. The above equality  implies that $X$ is connected. \qed

\begin{remark} The property  (c) in Proposition \ref{prop: tpp} is usually referred to as the  \emph{two-piece property} (TPP for  brevity).  The   subsets an Euclidean space satisfying this property are known as  \emph{$0$-tight sets}. \qed
\label{rem: tpp}
\end{remark}

Using the   above proposition  and Lemma \cite[Lemma 2.4]{Ku} we obtain the following result.

\begin{corollary} Suppose  that  $\Gamma$ is a   semialgebraic graph such that all its edges are straight  line segments.   Then  $\Gamma$ is tight if and only if there exists a  convex polyhedron $P$ such that   the following hold.

\smallskip

\noindent (a)  $\Gamma\subset P$.

\noindent (b)  The  vertices and edges of $P$ are contained in $\Gamma$.

\noindent (c) All the edges of $\Gamma$ are straight line segments.

\noindent (c) If $v$ is a vertex of  $\Gamma$ which is not a vertex of $P$ then  $v$ lies in the convex hull of its neighbors.\qed
\label{cor: s}
\end{corollary}

\begin{corollary} Any plane, closed, convex semi-algebraic curve  is tight.\qed
\end{corollary}

\section{Total curvature}
\setcounter{equation}{0}

Fix  a compact connected  semialgebraic graph $\Gamma\subset \bR^3$. For every integrable  function $S^2\ni\bu \mapsto f(\bu)\in \bR$ we
denote  by $\lan f(\bu)\ran\in \bR$ its average.  We would like to  investigate the average  $\lan w_\Gamma(\bu)\ran$ and  thus we would like to   know that the function
\[
S^2\setminus \Delta_\Gamma\ni \bu \mapsto w_\Gamma(\bu)\in \bZ
\]
is integrable.   This is  a consequence of the following result proved in the Appendix.

\begin{proposition} The function $S^2\setminus \Delta_\Gamma\ni\bu \mapsto w(\bu)\in \bZ$ is  semialgebraic. In particular, it is bounded and  measurable. \qed
\label{prop: morsesalg}
\end{proposition}

 Following \cite{CL} we define the \emph{total curvature} of  $\Gamma$ to be  average $K(\Gamma)$ of this function,
\[
K(\Gamma):=\lan w_\Gamma(\bu)\ran=\frac{1}{4\pi}\int_{S^2} |w(\bu)|d\si(\bu)|,
\]
where $|d\si|$ denotes the Euclidean area density on the unit sphere
$S^2$.  For the definition and (\ref{eq: morsein}) we obtain the
following generalization of  Fenchel's inequality.

\begin{corollary}
\[
K(\Gamma)\geq 1+b_1(\Gamma),
\]
with equality if and only if $\Gamma$ is tight.\qed
\label{cor: fench}
\end{corollary}

We  want  give a description of  $K(\Gamma)$ in terms of local geometric invariants  of $\Gamma$.  Observe that for any nondegenerate vector $\bu\in S^2$ we have
\[
w(\bu) =\# \Cr^1(\bu)+\sum_{v\in V}w(\bu,v)=\sum_{e\in E}\#\bigl(\,\Cr^1(\bu)\cap e\,\bigr)+\sum_{p\in V}w(\bu,p).
\]
  For every (open) edge  of $\Gamma$ we denote by $K(e)$ the total curvature
\[
K(e):= \frac{1}{\pi}\int_e|k_e(s)| |ds|,
\]
where $|ds|$ denotes the arclength  along $e$ and $k_e$ denotes the curvature along $e$.   Then (see \cite{CL} or \cite[\S 2]{Kui} ) we have
\[
\bigl\lan\, \# (\,\Cr^1(\bu)\cap e\,)\,\bigr\ran = K(e),
\]
so that,
\begin{equation}
\bigl\lan\, \# \Cr^1(\bu)\,\bigr\ran =\sum_{e\in E}K(e).
\label{eq: edges}
\end{equation}
For every vertex $p\in V$, and any nondegenerate unit vector $\bu$ we set
\[
\lambda(\bu, p)=\begin{cases}
2 & p\in V_{\min}(\bu)\\
d^-_p(\bu) &  p\in V\setminus V_{\min}(\bu).
\end{cases}
\]
Then
\[
w(\bu,p)= \lambda(\bu,p)-1\;\;\mbox{and}\;\;\sum_{p\in V} w(\bu,p)=\sum_{p\in V} \lambda(\bu, p) -\# V,
\]
so that,
\begin{equation}
\sum_{p\in V}\lan w(\bu,p)\ran= \sum_{p\in V}\lan \lambda(\bu,p)\ran -\# V.
\label{eq: avw}
\end{equation}

\begin{definition} For every point $p\in\Gamma$ we define
\[
\Sigma^+_\Gamma(p):=\bigl\{ \bu\in S^2;\;\;(\bu,\btau)>0,\;\;\forall  \btau\in \bsN_p\,\bigr\}.
\]
The closure of $\Sigma^+_\Gamma(p)$ is a  geodesic polygon on $S^2$, and we denote by $\si_\Gamma^+(p)$ its area.\qed
\label{def: angle}
\end{definition}

\begin{remark}Let us observe that if $p\in \Gamma$ is not a vertex, then  $\si_\Gamma^+(p)=0$.\qed
\label{rem: sangle}
\end{remark}

For every nondegenerate unit vector $\bu$,  and any vertex $p$ we have
\[
\lambda(\bu,p)+\lambda(-\bu,p)=\deg(p)+\begin{cases}
2 & \bu \in \Sigma^+_\Gamma(p)\cup-\Sigma^+_\Gamma(p)\\
0 & \mbox{otherwise},
\end{cases}
\]
and we conclude that
\[
2\lan \lambda(\bu,p)\ran=\lan \lambda(\bu,p)\ran +\lan \lambda(-\bu,p)\ran = \deg(p) +\frac{1}{\pi}\si^+_\Gamma(p).
\]
Hence,
\[
\sum_{p\in V}\lan \lambda(\bu,p)\ran =\frac{1}{2}\sum_{p\in V}\deg p +\frac{1}{2\pi}\sum_{p\in V}\si^+_\Gamma(p)=\# E+\frac{1}{\pi}\sum_{p\in V}\si^+_\Gamma(p).
\]
Using the equality (\ref{eq: avw}) we deduce
\[
\sum_{p\in V}\lan w(\bu,p)\ran =\# E-\# V+ \frac{1}{2\pi}\sum_{p\in V}\si^+_\Gamma(p)= \frac{1}{\pi}\sum_{p\in V}\si^+_\Gamma(p)-2+\bigl(1+b_1(\Gamma)\,\bigr).
\]
We have thus proved the following result.

\begin{theorem} The total curvature of the compact connected semialgebraic graph $\Gamma$ is given by
\[
K(\Gamma)=\frac{1}{2\pi}\sum_{p\in V} \si_\Gamma^+(p) + \sum_{e\in E}K(e)-\chi(\Gamma)=\frac{1}{2\pi}\sum_{p\in V} \si_\Gamma^+(p) + \sum_{e\in E}K(e)-2 + b_1(\Gamma)+1.
\]
In particular, the graph $\Gamma$ is tight if and only if
\[
\sum_{p\in V} \si_\Gamma^+(p) + 2\sum_{e\in E}\int_e |k_e(s)|\,|ds|=2\pi.\proofend
\]
\label{th: CL}
\end{theorem}

\begin{corollary} The total curvature  of $\Gamma$ is independent of the choice of vertex set $V$.\qed
\label{cor: cl}
\end{corollary}

\proof Indeed, Remark \ref{rem: sangle} implies that
\[
K(\Gamma)=\frac{1}{2\pi}\sum_{p\in \Gamma} \si_\Gamma^+(p) + \sum_{e\in E}K(e)-\chi(\Gamma).
\]
Neither one of the three summands   depends on the choice of vertex set.\qed

\begin{corollary}   Suppose $\Gamma$ is a compact connected semialgebraic  graph in $\bR^3$. For every  $\Gamma$-nondegenerate unit vector $\bu\in S^2$ we denote  by $\mu(\bu)$ the number of local minima of the function $h_{\bu}$ on  $\Gamma$, and by $\mu(\Gamma)$ the average, $\mu(\Gamma):= \lan \mu(\bu)\ran$.Then
\begin{equation}
\mu(\Gamma)= \frac{1}{2}\bigl(\, K(\Gamma)+\chi(\Gamma)\,\bigr)=\frac{1}{4\pi}\sum_{p\in V}\si^+_\Gamma(p)+\frac{1}{2\pi}\sum_{e\in E} \int_e |k(s)|\,|ds|.
\label{eq: av-min}
\end{equation}
\label{cor: min}
\end{corollary}

\proof Observe that
\[
2\mu(\bu)=2M_{\bu}(0)= M_{\bu}(1)+M_{\bu}(-1)= w(\bu)+\chi(\Gamma).
\]
The equality (\ref{eq: av-min})  is obtained by averaging the above identity. \qed

\begin{corollary} The total curvature of a plane, closed, convex semi-algebraic curve  is equal to $2$.\qed
\end{corollary}

\begin{corollary} Suppose $\Gamma$ is a planar, convex semialgebraic arc  with endpoints  $p_0\neq p_1$.  Let $\theta_0,\theta_1\in[0,\pi]$ denote the angles at $p_0$ and respectively $p_1$ between the arc $\Gamma$ and the line determined by $\pi_0$. Then
\[
\int_\Gamma|k(s)\,|ds|= \theta_0+\theta_1.
\]
\end{corollary}

\proof  Denote by $\hat{\Gamma}$ the closed curve obtained from $\Gamma$ by  connecting $p_0$ to $p_1$. Then $\hat{\Gamma}$ is  a plane convex curve and using Theorem \ref{th: CL} we deduce
\[
2=K(\hat{\Gamma})=\frac{1}{\pi}\int_\Gamma|k(s)\,|ds|+ \frac{1}{2\pi}\bigl(\, (2\pi-2\theta_0)+(2\pi-2\theta_1)\,\bigr)
\]
which implies the claimed equality.\qed

\begin{corollary}  Suppose the graph $\Gamma$ is contained in a plane and it  is the union of a closed, convex semialgebraic curve $B$ and a finite union  of line segments   contained in the  region $R$ bounded by $B$ such  that  $R\setminus \Gamma$ is a finite union of convex sets. Then $\Gamma$ is tight.
\label{cor: c}
\end{corollary}

\proof  The curve $B$ is a disjoint union of  (open) edges and vertices of $B$.   Denote by $V_B$ (respectively $E_B$) the collection of vertices (respectively  open) edges contained in $B$ and by $V'_B$ (respectively $E'_B$)  the collection of vertices (respectively open edges)  not contained in $B$.  If $v\in V'_B$  then the cone generated by $\bsN_v$ is the plane of the graph $\Gamma$ so that $\si_\Gamma^+(v)=0$. Similarly, if $e\in E'_B$, then $K(e)=0$ since $e$ is a straight line segment. Hence
\[
K(\Gamma)=\frac{1}{2\pi}\sum_{v\in V_B}\si^+_\Gamma(B) +\sum_{e\in E_B} K(e)-2\chi(\Gamma)=K(B)-\chi(\Gamma)=2-\chi(\Gamma)=1+b_1(\Gamma).
\]
This proves that $\Gamma$ is tight. \qed

\begin{remark} Observe that when $\Gamma$ is a polygonal simple close curve then the  formula (\ref{eq: av-min}) specializes to  Milnor's formula \cite[Thm. 3.1]{Mil}.

(b) Recently,  Gulliver and  Yamada  have proposed in \cite{GY} a different notion of total curvature.   For every nondegenerate vector $\bu$  and every vertex $p$ they define a defect at  $p$ to be the integer
\[
\delta(\bu, p)= \bigr(\,d_p^-(\bu,p)-d_p^+(\bu)\,\bigl)^+,\;\;x^+ :=\max(x,0),\;\;.
\]
The Gulliver-Yamada total curvature is then the real number $T(\Gamma)$ defined by
\[
\frac{1}{\pi}T(\Gamma)=\sum_{p\in V}\lan\, \delta_+(\bu,p)\,\ran +\sum_{e\in E}K(e).
\]
Observe that
\[
\delta(\bu,p)+\delta(-\bu,p)=|d_p^-(\bu,p)-d_p^+(\bu)|
\]
so that
\[
\frac{1}{\pi}T(\Gamma)= \frac{1}{2}\sum_{p\in V}\bigl\lan\, |d_p^-(\bu,p)-d_p^+(\bu)|\,\bigr\ran +\sum_{e\in E}K(e)
\]
Let us show that when  all vertices have degrees $\leq 3$, then
\begin{equation}
\frac{1}{\pi}T(\Gamma)=K(\Gamma).
\label{eq: GY}
\end{equation}

Indeed, for  every vertex $p\in V$, and every  unit vector $\bu$, such that both $\bu$ and $-\bu$ are regular,  we have
\[
\delta(\bu,p)+\delta(-\bu,p)=w(p,\bu)+ w(p,-\bu)=\begin{cases}
3 & \bu\in \Sigma^+_p(\Gamma)\cup -\Sigma_p^+(\Gamma)\\
1 & \mbox{otherwise}.
\end{cases}
\]
The equality (\ref{eq: GY}) follows by averaging the above identity.

In general $K(\Gamma)\neq \frac{1}{\pi}T(\Gamma)$. To see this  consider the  graph $\Gamma_0$ obtained  by joining the north pole to the south pole of the unit sphere by   $n$  meridians  of longitudes $\frac{2k\pi}{n}$, $k=1,\dotsc, n$. If $n$ is an odd  integer then
\[
|d_p^-(\bu,p)-d_p^+(\bu)|=1
\]
for any vertex $p$ and almost all $\bu$'s.   In this case we have $\si_p^+(\Gamma_0)=0$ for any vertex $p$ and we deduce
\[
K(\Gamma_0)= \sum_{e\in E}K(e) +n-2.
\]
On the other hand
\[
\frac{1}{\pi}T(\Gamma_0)= 1+ \sum_{e\in E}K(e).
\]

(c) Taniyama \cite{Tan}  has proposed another notion of total curvature $T_1(\Gamma)$ given by
\[
\frac{1}{\pi}{T_1(\Gamma)}=\sum_{e\in E}K(e)+\sum_{p\in V} \theta(p),
\]
where for every vertex $p$ the  quantity $\theta(p)$ is the sum
\[
\sum_{\btau_0\neq\btau_1\in \bsN_p} \bigl(\,\pi -\measuredangle(\btau_0,\btau_1)\,\bigr).
\]
We can check by direct computation that $K(\Gamma_0)\neq \frac{1}{\pi}T(\Gamma_0)\neq \frac{1}{\pi}T(\Gamma_1)$.   It seems that  Taniyama's definition is closer in the spirit to the approach of I. F\'{a}ry \cite{Fa} where he showed that the total curvature  of  simple closed curve in $\bR^3$ is equal to the average of the total curvatures of its  projections on all the two dimensional planes. \qed
\label{rem: curv}
\end{remark}

From  (\ref{eq: GY}) and  Corollary \ref{cor: fench}  we obtain the following generalization of the first half of \cite[Thm.2]{GY}.

\begin{corollary} Suppose $\Gamma$ is a connected,  semialgebraic graph  such that each of its  vertices  has degree $\leq 3$. If $T(\Gamma)$ denotes the   Gulliver-Yamada total curvature of $\Gamma$ then
\[
T(\Gamma)=\pi K(\Gamma)\geq \pi\bigr(\,1+b_1(\Gamma)\,\bigl).\proofend
\]
\label{cor: GY}
\end{corollary}

\section{Tightness}
\setcounter{equation}{0}
We want to give a complete explicit classification of tight semi-algebraic  graphs in $\bR^3$.

\begin{theorem} A semialgebraic graph $\Gamma$ is tight if an only if  its of one the following types: type C described in Corollary \ref{cor: c},  and  type S described in Corollary \ref{cor: s}.
\label{th: tight}
\end{theorem}

\proof    Suppose  $\Gamma$ is a connected semialgebraic graph. As usual we denote by $V$ the set of vertices and by $E$ the set of edges.  We denote by $\bsR_\Gamma\subset S^2$ the set of $\Gamma$-regular unit vectors.

 Assume $\Gamma$ is tight. According to Proposition \ref{prop: tpp} this means that the intersection of $\Gamma$ with any closed half-space is connected.

For every point $p$ on an (open) edge  $e$ we denote by $L_p$ the affine  line tangent to the edge  $e$ at $p$ and   by $k(p)$ the curvature of $e$ at $p$. If $k(p)\neq 0$ we denote $\bn(p)$ the normal vector to $e$ at $p$. We set
\[
\eC(\Gamma):=\bigl\{ p\in \Gamma\setminus V;\;\; k(p)\neq 0\}. \bigl(\, V\cup F(\Gamma)\,\bigr),\;\;\eC(e):=\eC(\Gamma)\cap e,\;\;\forall e\in E.
\]
We will refer to $\eC(\Gamma)$ as the \emph{curved region} of $\Gamma$ since it consists of the points along edges where the curvature is nonzero.  The set  $\eC(e)$ is semialgebraic,  open in $e$ and consists of finitely many arcs.

If $\eC(\Gamma)=\emptyset$ then all the edges of $\Gamma$ are straight line segments and   the theorem reduces to  Corollary \ref{cor: s}. Thus  in the sequel  we will assume that $\eC(\Gamma)\neq \emptyset$.

\begin{lemma} Suppose $p$ is a point  on an open edge of the tight graph $\Gamma$ such that $k(p)\neq 0$. Then the graph $\Gamma$ is contained  in the affine osculator plane, i.e., the affine plane  determined by the affine tangent line $L_p$ and the normal vector $\bn(p)$, and it is situated entirely in one of the closed half-planes determined by $L_p$.
\label{lemma: curv}
\end{lemma}

\proof  Consider  a unit vector $\bu$ such that $(\bu,\bn(p))>0$.  Set $c=h_{\bu}(p)$  The intersection between the half-space $\{h_{\bu}\leq c\}$ and $\Gamma$ is closed, connected and semialgebraic. The curve selection theorem implies that  and $p$ is an isolated point of this intersection.  Hence the intersection consists of  a single point so that $\Gamma$ is contained in the half-space
\[
H^+_{\bu,p}=\bigl\{ q\in \bR^3;\;\;(\bu,q-p)\geq 0\}.
\]
Thus
\[
\Gamma \subset  \bigcap_{(\bu,\bn(p)>0} H^+_{\bu,p}.
\]
This intersection is  a half-plane in the osculator plane determined by the affine tangent line $L_p$. \qed

Using  a similar argument  as in the proof  of Lemma \ref{lemma: curv}  we deduce the following result.

\begin{lemma} The edges  of $\Gamma$ are convex arcs.\qed
\label{lema: conv}
\end{lemma}

For every vertex $v$ of $\Gamma$ we denote by $\eC_v$ the convex cone spanned by the vectors in $\bsN_v$ and by $\eC_v^*$ is dual cone
\[
\eC_v^*=\bigl\{  \bp\in \bR^3;\;\; (\bp,\btau)\geq 0;\;\;\forall \btau\in \bsN_v\,\}.
\]
The intersection of $\eC_v^*$ with the unit sphere is the closure of the open region $\Sigma^+_\Gamma(v)$ introduced in Definition \ref{def: angle}.

\begin{lemma} For every vertex $v$ of $\Gamma$ we have $\Gamma\subset \eC_v$.
\label{lemma: angle}
\end{lemma}

\proof  Note first that $\Sigma^+_\Gamma(v) =\emptyset$ if and only if $\si^+_\Gamma(v)=0$.   If $\si_\Gamma^+(v)=0$ then $\eC_v=\bR^3$ and the statement is obvious.  Assume $\si_\Gamma^+(v)>0$.

Then  $\Sigma^+_\Gamma(v)$ is an open subset   of $S^2$ and $\Sigma^+_\Gamma(v)\setminus \Delta_\Gamma$ is dense in $\Sigma^+_\Gamma(v)$.  For any  $\bu\in \Sigma^+_\Gamma(v)\setminus \Delta_\Gamma$  the function $h_{\bu}$ has a local  minimum at $v$, and since $\Gamma$ is tight, it has to be  the  absolute minimum. In other words $\Gamma$ is contained  in the  half-space $H^+_v(\bu)$ so that
\[
\Gamma\subset \bigcap_{\bu\in \Sigma^+_\Gamma(v)\setminus \Delta_\Gamma}H^+_v(\bu)= \bigcap_{\bu\in \Sigma^+_\Gamma(v)} H^+_v(\bu)
\]
(cl=closure)
\[
=\bigcap_{\bu\in {\rm cl}\,(\Sigma^+_\Gamma(v))} H^+_v(\bu)=(\eC^*_v)^*=\eC_v,
\]
where at the last step we have used the fact that a closed convex cone coincides with  its bidual. \qed

We set
\[
V_{ext}=V_{ext}(\Gamma):=\bigl\{v\in V;\;\;\si^+_\Gamma(v)>0\,\bigr\},
\]
and we will refer to the vertices in $V_{ext}$ as \emph{extremal vertices}.

\begin{lemma} Suppose $\eC(\Gamma)\neq\emptyset$.  Then there exists a  semialgebraic closed convex curve  $B$ such that the following hold.

\smallskip

\noindent (a) $B\subset \Gamma$.

\noindent (b) $\Gamma$ is contained in the convex hull $R$  of $B$.

\noindent (c) $\Gamma\setminus B$ is a union of straight line segments.

\noindent (d) The components of $R\setminus \Gamma$ are convex open sets.

\label{lemma: plane}
\end{lemma}

\begin{figure}[h]
\centerline{\epsfig{figure=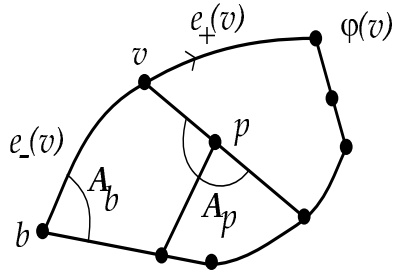,  height=1.4in, width=2in}}
\caption{\sl A tight graph and its set $V'_{ext}$.}
\label{fig: 2}
\end{figure}

\proof     Since $\eC(\Gamma)\neq \emptyset$  Lemma \ref{lemma: curv} implies  that $\Gamma$ must be a planar.  For every $v\in V$ we denote by $\bsA_v$ the intersection of $\eC_v$  the unit circle in  the plane containing $\Gamma$ of the curve.  Set $\theta_v:={\rm length}\,(\bsA_v)$.  Denote  by $V'_{ext}$ the set  of vertices  such that  $\theta_v\leq \pi$. Clearly $V_{ext}\subset V'_{ext}$ but in general the inclusion is strict.  For example, the vertices $p$ and $v$ in Figure \ref{fig: 2} belong to $V'_{ext}$, but not to $V_{ext}$. The vertex $b$ belongs to $V_{ext}$.

Note that if $\theta_v>\pi$ then $\theta_v=2\pi$. Moreover
\[
\si^+_\Gamma(v)=  2(\pi-\theta_v).
\]
 For $v\in V'_{ext}$ the angle $\bsA_v$ is spanned by two vectors  $\btau_\pm\in \bsN_v$,   where $\btau_\pm$ are ordered such that  the counterclockwise angle from $\btau_-$ to $\btau_+$  is $\leq\pi$.   The two  vectors $\btau_\pm$ correspond to two edges $e_\pm(v)$ incident to $v$; see Figure \ref{fig: 2}.

If we start at $v$ and travel along the edge $e_+(v)$ we   encounter another vertex $\vfi(v)$ of $\Gamma$.  From Lemma \ref{lemma: curv} and \ref{lemma: angle} we deduce that during our travel the graph  $\Gamma$ is situated to the right of the  affine tangent lines to $e_+(v)$ oriented by the  direction of motion; see Figure \ref{fig: 2}.  This implies  that $\vfi(v)\in V'_{ext}$ and $e_+(v)=e_-(\vfi(v))$.

We have thus obtained a  map $\vfi: V'_{ext}\ra V'_{ext}$ such that $v$ is connected to $\vfi(v)$ by the edge $e_+(v)$ and  $\Gamma$ is situated to the right of the edge $e_+(v)$  oriented  by the motion from $v$ to $\vfi(v)$.

Suppose $S\subset V_{ext}'$ is a minimal $\vfi$-invariant subset  of $V_{ext}$. Then  we can label the vertices in $S$ as $s_1,s_2,\dotsc, s_k$, so that
\[
k=|S|,\;\;v_2=\vfi(v_1),\dotsc,\;\;v_k=\vfi(v_{k-1}),\;\;v_1=\vfi(v_k).
\]
 The succession of edges $e_+(v_1),\dotsc, e_+(v_k)$ determines a closed, clockwise oriented curve. It is convex because    it is situated on one side    of the affine tangent lines to the smooth points of this curve and it is contained in each of the angles $\bsA_{v_i}$.  Denote by $B$ this closed convex curve.

 The graph $\Gamma$ is contained in the region $R$ bounded by  $B$.  If $e$ is an (open) edge of $\Gamma$ not intersecting $B$, then $e$ must be a  line segment. Indeed, if $p\in e$ is a point where $k(p)\neq 0$, then the line $\eN_p$ through $p$ determined by $\bn(p)$ intersects   $B$ in two different points $q_1,q_2$ which, according to Lemma \ref{lemma: curv}, are situated in  the  same half-plane determined by the affine tangent line $L_p$. In particular, the point $p$ is not contained  on the closed segment $[p_1,p_2]$. This is a contradiction  because  the intersection of the line $\eN_p$ with $R$ is precisely the segment $[p_1,p_2]$ which implies that $\Gamma$  has a point $p$ not contained in $R$. Arguing in a similar fashion we deduce that $V_{ext}\subset B$.

 For every vertex $v\in V$ the  set $\bsN_v$ is contained in the unit circle in the plane of $\Gamma$. Thus, we can cyclically counterclockwise  order the vectors in $\bsN_v$. If $v\in V\setminus V'_{ext}$ then the (counterclockwise) angle between two consecutive vectors in $\bsN_v$ (with respect to this cyclic ordering) is $<\pi$ because   in this case  the planar cone spanned by $\bsN_v$  must coincide with the plane of $\Gamma$.

 Suppose $D$ is a  component of $R\setminus \Gamma$. Its boundary is a disjoint union of   finitely many  vertices and open edges of $\Gamma$.   We fix the clockwise orientation of $\pa D$.  At each vertex $v\in V\cap \pa D$ we have two edges an incoming and an outgoing edge contained in $\pa D$.    They determine two vectors $\btau_0,\btau_1\in \bsN_v$.  Since $D$ is a connected component of the complement of $\Gamma$ in $R$, these two vectors must be successive vectors in $\bsN_v$ with respect to the  counterclockwise cyclic order on $\bsN_v$.

 If $v\in B\cup V'_{ext} $ then the angle between these vectors is $\leq \pi$ because the curve $B$ is convex. If $v$ is in the interior of $R$ and  $v\in V\setminus V'_{ext}$ then  as we have seen above,   the counterclockwise angle between these vector must be smaller than $\pi$. This proves that  $D$ is a convex open region.\qed

This completes the proof of Theorem \ref{th: tight}. \qed

\begin{remark} Because  tight sets satisfy  the two-piece-property (see Remark \ref{rem: tpp})  we can use  the characterization   of planar $0$-tight planar sets  in   \cite[Thm. 1.3]{Kui} to give an alternate proof  of Lemma \ref{lemma: plane}.\qed
\end{remark}

\begin{remark} We chose to work in the category of semialgebraic sets because they may be familiar to a larger audience.  In fact  all the results in this paper are valid in any $o$-minimal category.  We refer to \cite{Dr1} for more information about $o$-minimal sets and maps.\qed
\end{remark}

\section{Knottedness}
\setcounter{equation}{0}
As we mentioned in the introduction,    Milnor proved that  the
knotting of a circle $C\subset \bR^3$ requires  a substantial
amount  of curvature. More precisely,  if $\mu(C)<2$ then $C$
cannot be knotted. Gulliver  and Yamada  \cite{GY} extended this
result to singular situations.  Let us define a $\theta$-graph to
be a graph homeomorphic to the union of a round circle and a
diameter.  Such a graph is $3$-regular, and  for these graphs  the
Gulliver-Yamada total curvature coincides with the notion of total
curvature  introduced in this paper.   Theorem 2 in \cite{GY} can
be rephrased as follows. If $\Gamma\subset \bR^3$ is  a
semialgebraic graph homeomorphic to a $\theta$-graph, and
$\mu(\Gamma)<\frac{3}{2}$, then the embedding of $\Gamma$ is
isotopic to a planar embedding.

\begin{figure}[h]
\centerline{\epsfig{figure=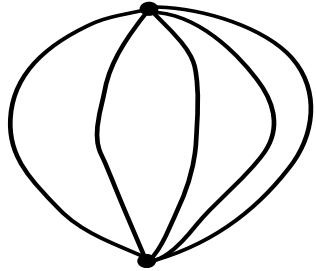,  height=1.4in, width=1.6in}}
\caption{\sl The suspension of $5$ points.} \label{fig: 3}
\end{figure}

Denote by $\Sigma_n$ the suspension  of an $n$-point set. In
Figure \ref{fig: 3} we have  depicted a planar embedding of
$\Sigma_5$.   We will refer to the two distinguished points of the
suspension as the \emph{poles}.   Observe that $\Sigma_2$ is a
circle, and $\Sigma_3$ is a $\theta$-graph. Our next result
generalizes  the unknottedness results of  Milnor \cite{Mil} and
Gulliver-Yamada \cite{GY}.

\begin{theorem} Suppose $\Gamma\subset \bR^3$ is a semialgebraic  subset   homeomorphic to the suspension of $n$-points.  If $\mu(\Gamma) <2$ then $\Gamma$ is unknotted, i.e., it is isotopic  with a planar  embedding of the suspension.
\label{th: unknot}
\end{theorem}

\proof     We follow the strategy pioneered by J. Milnor in
\cite{Mil}.  Fix a  vertex set $V$ on $\Gamma$. The poles belong
to the vertex set.    The complement of the poles in $\Gamma$ is a
semialgebraic set with $n$ connected components. We will refer to
these  as the      \emph{meridians} of the suspension.  They are
(open) semialgebraic  arcs without self intersections.

Denote by $\eU\subset S^2$ the set consisting of
$\Gamma$-nondegenerate unit vectors $\bu$ such that   $h_{\bu}$
has at least two local minima on $\Gamma$.  The set  $\eU$ is
semi-algebraic.  Since $\mu(\Gamma) <2$ we deduce  that
\[
{\rm Area}\,(\eU)<\frac{1}{2}{\rm Area}\,(S^2).
\]
If $\eM:=S^2\setminus \eU$ then $\eM$ is semialgebraic and ${\rm
Area}\,(\eM)>\frac{1}{2}{\rm Area}\,(S^2)$.    Hence
\[
{\rm Area}\,\bigl(\,\eM\cap (-\eM)\,\bigr)>0.
\]
Since $\eM\cap (-\eM)$ is  semialgebraic, and has nonzero area, we
deduce that it has nonempty interior.  We denote this interior by
$\eV$.  For any $\bu\in \eV$ the function  $h_{\bu}$   has a
unique local minimum and a unique  local maximum  on  $\Gamma$.
From the results in \cite{Pig} we deduce that   we   can choose
$\bu\in \eV$ satisfying the additional condition that the
restriction of $h_{\bu}$ on $\Cr(\bu)$ is injective.

The poles of $\Gamma$ are  critical points of $h_{\bu}$ so that
$h_{\bu}$ has distinct values on these poles. We  label the poles
by $P_\pm$ so that $h_{\bu}(P_+)>h_{\bu}(P_-)$. We set $m_\pm
:=h_{\bu}(P_\pm)$. We need to distinguish two cases.

\smallskip

\begin{figure}[h]
\centerline{\epsfig{figure=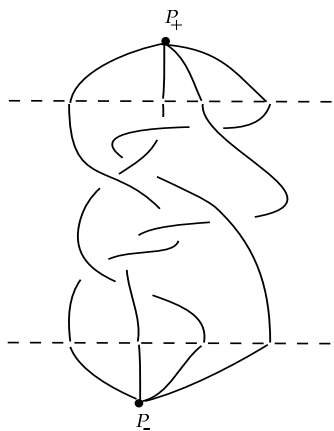,  height=2.2in, width=1.7in}}
\caption{\sl A braided embedding of $\Sigma_4$.} \label{fig: 8}
\end{figure}

\noindent {\bf A.}  \emph{ The poles are the only local extrema of
$h_{\bu}$ on $\Gamma$.}  In this case, the restriction  of
$h_{\bu}$  on each meridian   induces a semialgebraic homeomorphic
between that meridian and the open interval $(m_-,m_+)$.  The
graph is then  a braid with the top and bottom capped-off; see
Figure \ref{fig: 8}.

\begin{figure}[h]
\centerline{\epsfig{figure=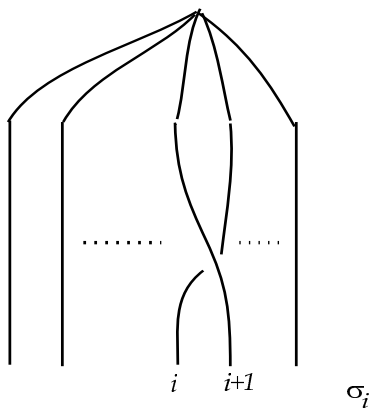,  height=2.1in, width=1.9in}}
\caption{\sl The elementary braid $\si_i$ with the top capped
off.} \label{fig: 9}
\end{figure}

Such a braided embedding is isotopic to a planar  embedding
because a braid with the top  capped off can be untwisted by an
isotopy   which keeps the bottom of the braid fixed. This can be
see easily for the   basic braids $\si_i^\pm$; see Figure
\ref{fig: 9}.   By Artin's classical result (see e.g. \cite[Thm.
1.8]{Bir})  any braid is a composition of such elementary
braids. Thus,   we can  inductively untwist any capped braid.

\noindent {\bf B.}  \emph{One of the poles $P_\pm$ is not an
absolute extremum of $h_{\bu}$}.    We will reduce this case to
the previous one. More precisely, we will prove that we can isotop
$\Gamma$   an embedding  of $\Sigma_n$ such that $h_{\bu}$  is a
stratified Morse function with a unique relative maximum at $P_+$
and a unique relative minimum located at $p_-$.

\begin{figure}[h]
\centerline{\epsfig{figure=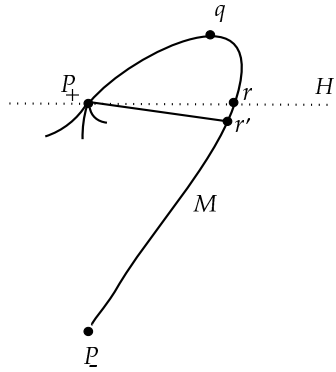,  height=1.9in, width=1.7in}}
\caption{\sl Deforming $\Gamma$ to a braided embedding.}
\label{fig: 10}
\end{figure}

Suppose that the unique relative maximum is achieved at a point
$q$ located on one of the meridians. Denote this meridian by $M$.
The hyperplane $H$  determined by $\bu$ and containing $P_+$
intersects the meridian $M$ at a unique point $r$; see Figure
\ref{fig: 10}.  (If there were several points of intersection with
$M$,  the function $h_{\bu}$ would have several local maxima on
this meridian.)

The pole $P_+$ is a stratified critical point and, since the
function $h_{\bu}$ does not have  local maxima on  meridians
different from $M$ we deduce that  we can find an open, convex
polyhedral cone $\eC$  with vertex at $P_+$   satisfying the
following  conditions.

\begin{itemize}

\item The cone $\eC$ is situated  below the hyperplane $H$, i.e.,
it is contained in the open half-space $\{h_{\bu}<m_+\}$.

\item The cone $\eC$  contains all the meridians other than $M$.

\end{itemize}

We deduce that if  a point $r'$ on the meridian $M$ is
sufficiently close to $r$, then the open line segment $(P_+r')$
does not intersect any of the meridians situated below $H$. Now
choose $r'\in M$ sufficiently close to $r$ and situated slightly
below $H$; see Figure \ref{fig: 10}. Moreover, we can assume that
$M$ is smooth at $r'$ and  thus $r'$ is not a critical point of
$h_{\bu}$. Clearly we can deform the arc $P_+qr'$ to the segment
$[P_+r']$ while keeping the endpoint fixed such that  during the
deformation we do not intersect the meridians different from  $M$.
We   have thus isotoped $\Gamma$ to an embedding with the property
that $h_{\bu}$  is a stratified Morse function with  no local
maxima along the meridians.

Arguing similarly, we can isotop $\Gamma$ such that $h_{\bu}$ has
no local minima along the meridians.  Thus we  have isotoped
ourself to case  {\bf A}. This  concludes the proof of Theorem
\ref{th: unknot}. \qed

 \begin{figure}[h]
\centerline{\epsfig{figure=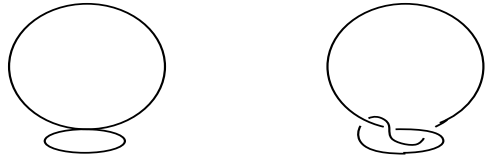,  height=0.8in, width=2.5in}}
\caption{\sl A surgery  from the singular figure eight to a trefoil knot.} \label{fig: 11}
\end{figure}

\begin{remark}(a)  The special case $b=2$ is due to F\'{a}ry \cite{Fa} and  Milnor \cite{Mil}. The case  $n=3$ was investigated by  Gulliver and Yamada \cite{GY}  where they established the unknotedd ness of a theta graph  whose crookedness satisfies $\mu<\frac{3}{2}$.

\noindent (b)   The upper bound  $2$  on crookedness  in Theorem
\ref{th: unknot} cannot be improved.  For example, if $\Gamma$ is
the figure eight curve in the left-hand side    of Figure
\ref{fig: 11} then $\mu(\Gamma)=2$.  We can  find surgeries of
this   figure eight curve producing the trefoil knot in the
right-hand side while barely changing    the total curvature so
that $\mu({\rm trefoil})\approx 2$.    Essentially, this surgery
replaces  a   straight line segment with an arc of circular helix
with the same endpoints, and axis parallel to the segment. If the
slope of  of this arc of helix is very large then  the total
curvature is very small.  More precisely the total  curvature of
the arc
\[
[0, 2\pi c]\ni s \mapsto \bigl(\,a\cos\frac{s}{c}, a\sin \frac{s}{c},\frac{bs}{c}\,\bigr)\in \bR^3,\;\; c=\sqrt{a^2+b^2},
\]
is $\frac{2a}{c}=\frac{2}{\sqrt{1+m^2}}$ where $m=\frac{b}{a}$ is  its slope. \qed \end{remark}

\begin{remark} The above result  is not valid for   graphs that are not
suspensions of finite sets. In fact the situation is dramatically
different. For such   graphs, it is possible that they are almost
tight $\mu(\Gamma)$ is very close to $1$ and still be knotted. We
describe below such an instance.

 \begin{figure}[h]
\centerline{\epsfig{figure=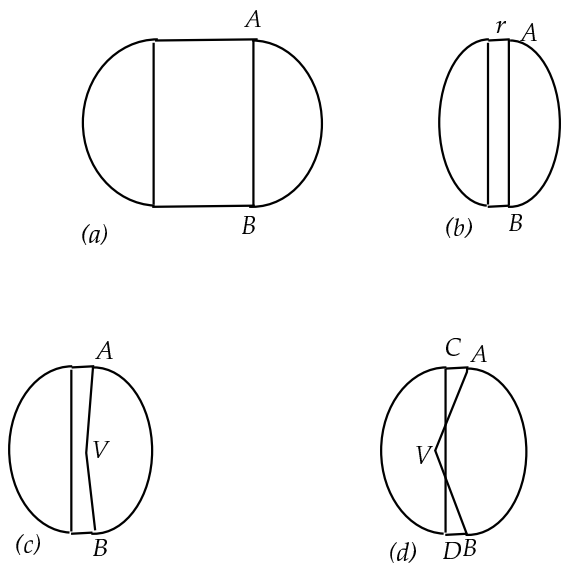,  height=2.9in, width=2.9in}}
\caption{\sl Deforming a simple tight set.} \label{fig: 5}
\end{figure}

Consider the planar  semialgebraic set $\Gamma_a$ depicted  in
Figure \ref{fig: 5}(a).  In this case  $\mu(\Gamma)=1$.  We shrink
the horizontal edges  to obtain a set  $\Gamma_b$ as  in Figure
\ref{fig: 5}(b)  still satisfying $\mu(\Gamma_b)=1$. We denote by
$r$ the common length of the two horizontal edges.

\begin{figure}[h]
\centerline{\epsfig{figure=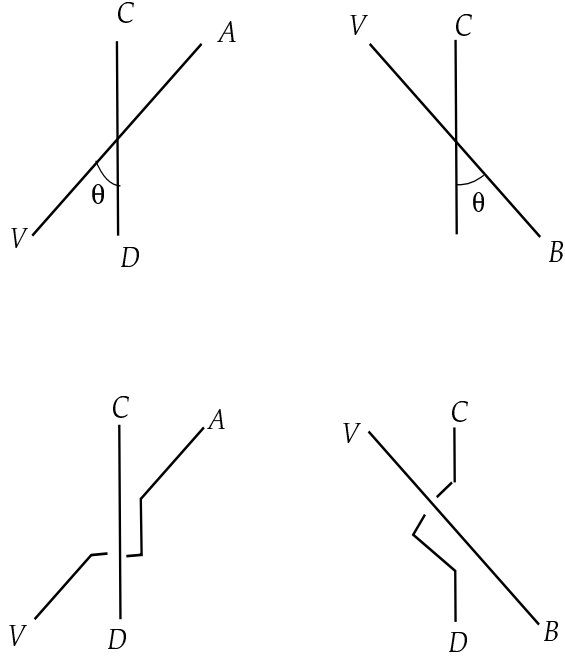,  height=3.4in, width=2.9in}}
\caption{\sl Resolving a crossing as an under/overcrossing.}
\label{fig: 6}
\end{figure}

Now deform the vertical edge $AB$ to obtain the broken   line
$AVB$ as in Figure \ref{fig: 5}(c). The triangle $AVB$ is
isosceles.  The new semialgebraic set $\Gamma_c$ satisfies
\[
\mu(\Gamma_c)= 1+\frac{\theta}{\pi},
\]
where $\theta$ denotes the angle between $AB$ and $VA$.  Continue
deforming the   broken line  $AVB$    until $V$  crosses the
vertical chord $CD$, and push $V$ a bit more to the  left to
obtain the semialgebraic set  $\Gamma_d$ depicted in Figure
\ref{fig: 5}(d).

If we continue to denote by $\theta$ the (new) angle between $VA$
and $AB$ then
\[
\mu(\Gamma_d)=1+\frac{\theta}{\pi}.
\]
On the left-had side of Figure \ref{fig: 6} we describe how to
replace the intersection between  the lines $VA$ and $CD$  with by
creating a small  (under) ``dimple'' in  the arc $AB$.  The new
arc $VA$ will  go below $CD$. We can arrange that   resulting
change in $\mu(\Gamma)$ due to this dimple  is $<\frac{\ve}{4}$.
In the right-hand side  of Figure \ref{fig: 7} we  replace
similarly the intersection of $VB$ and $CD$ with a crossing of
$VB$ over the new arc $CD$. Again we can arrange that the change
in $\mu(\Gamma)$ due to this surgery is $<\frac{\ve}{4}$.

After all these transformations we obtain  a subset in $\bR^3$
homeomorphic to the original  set $\Gamma_a$ and isotopic to the
set depicted in Figure \ref{fig: 7}. We observe two disjoint
circles  forming a Hopf link.

\begin{figure}[h]
\centerline{\epsfig{figure=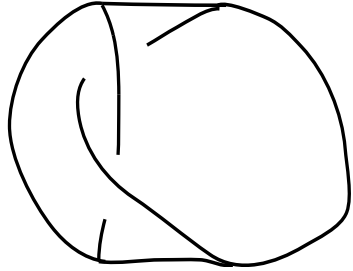,  height=1.4in, width=1.8in}}
\caption{\sl The formation of a Hopf link} \label{fig: 7}
\end{figure}

By choosing   the length $r$ in Figure \ref{fig: 5}(b)
sufficiently small, we can arrange that the angle $\theta$ between
$VA$ and $CD$ in Figure \ref{fig: 5}(d) is $<\frac{\pi\ve}{2}$.
After performing all these operations we obtain    a semialgebraic
$\Gamma$ set   isotopic with the set in Figure \ref{fig: 7}  and
satisfying
\[
\mu(\Gamma)
<1+\frac{\theta}{\pi}+\frac{\ve}{4}+\frac{\ve}{4}<1+\ve.
\]
This set cannot be unknotted to the set $\Gamma_a$ because of the
presence of the nontrivial Hopf link.\qed
\label{rem: knot}
\end{remark}


\appendix

\section{Basic real algebraic geometry}
\setcounter{equation}{0}
We   want to present a few basic facts of real algebraic geometry
used throughout the paper. For proofs and more details we refer to
\cite{BCR, Dr1,Co}.

A subset  $A\subset \bR^n$ is called \emph{semialgebraic} if it is
a finite union of sets  described by finitely many polynomial
inequalities. We denote by $\eS^n$ the collection of
semi-algebraic subsets of $\bR^n$ and we set $\eS=\cup_{n\geq
1}\eS^n$.

Using the canonical embedding $\bR^n\subset \bR^{n+1}$ we can
regard $\eS^n$ as a subcollection  of $\eS^{n+1}$.  If $A\in
\eS^m$ and $B\in \eS^n$, then a (possibly discontinuous) map
$f:A\ra B$ is   semialgebraic if its graph $\Gamma_f\subset
A\times B$ is a semialgebraic set.  Sometimes we will refer to
the semialgebraic sets/functions as  \emph{definable} sets/
functions.

We list below  some of the basic properties of the semialgebraic
sets.

\noindent $\bullet$  The collection $\eS$ is closed under boolean
operations: union,  intersection complement, cartesian product.

\noindent $\bullet$ (Tarski-Seidenberg) If $T:\bR^m\ra \bR^n$ is
an affine map, then $T(\eS^m)\subset \eS^n$.

\noindent $\bullet$ The  image and preimage of a definable set via
a definable map is a definable set.

 \noindent $\bullet$ (\emph{Piecewise  smoothness of  one variable semialgebraic functions.}) If $f:(0,1)\ra \bR$ is a semialgebraic function, then there exists
 \[
0=a_0< a_1<a_2<\cdots <a_n=1
\]
such that the restriction of  $f$ to each subinterval
$(a_{i-1},a_i)$ is real analytic and   monotone. Moreover $f$
admits right and left limits at any $t\in[0,1]$.

\noindent $\bullet$ (\emph{Closed graph theorem.})  Suppose $X$ is
a  semialgebraic set and $f: X\ra \bR^n$ is a  semialgebraic
bounded function.  Then $f$ is continuous if and only if  its
graph is closed in $X\times \bR^n$.

\noindent $\bullet$ (\emph{Curve selection.}) If $A$ is a
definable  set, and $x\in\cl(A)\setminus A$, then there exists an
$\eS$ definable  continuous map
\[
\gamma:(0,1)\ra A
\]
such that $x=\lim_{t\ra 0} \gamma(t)$.

\noindent $\bullet$  Any definable set has finitely many
connected components, and each of them  is definable.

\noindent $\bullet$   Suppose $A$ is a definable  set, $p$ is a
positive integer, and $f: A\ra \bR$ is a definable function. Then
$A$ can be partitioned into finitely many  definable sets
$S_1,\dotsc, S_k$,     such that each  $S_i$ is a $C^p$-manifold,
and each of the restrictions $f|_{S_i}$ is a $C^p$-function.

\noindent $\bullet$ (\emph{Triangulability.})  For every   compact
definable set $A$, and any finite collection of definable  subsets
$\{S_1,\dotsc, S_k\}$, there exists  a compact simplicial complex
$K$, and a  definable homeomorphism
\[
\Phi: K\ra A
\]
such that  all the sets $\Phi^{-1}(S_i)$ are unions of  relative
interiors of faces of $K$.

\noindent $\bullet$ (\emph{Definable selection}.)  Suppose
$A,\Lambda$ are definable. Then a \emph{definable} family of
subsets of $A$ parameterized by $\Lambda$ is a definable subset
\[
S\subset A\times \Lambda.
\]
We set
\[
S_\lambda:= \bigl\{ a\in A;\;\; (a,\lambda)\in S\,\bigr\},
\]
and we denote by $\Lambda_S$ the projection of $S$ on $\Lambda$.
Then there exists a definable function $s:\Lambda_S\ra A$ such
that
\[
s(\lambda)\in S_\lambda,\;\;\forall \lambda\in \Lambda_S.
\]

\noindent $\bullet$ (\emph{Dimension.})  The  dimension of a
definable  set $A\subset \bR^n$ is the supremum over all the
nonnegative integers $d$ such that there exists a $C^1$
submanifold of $\bR^n$ of dimension $d$ contained in $A$.  Then
$\dim A <\infty$,  $\dim (X\times Y)=(\dim X)(\dim Y)$, $\forall
X,Y\in \eS$, and
\[
\dim (\cl(A)\setminus A) <\dim A.
\]
Moreover, if $(S_\lambda)_{\lambda\in\Lambda}$ is a definable
family of definable sets then the function
\[
\Lambda\ni\lambda \mapsto\dim S_\lambda
\]
is definable.

\noindent $\bullet$   If $f: A\ra B$ is a semialgebraic
bijection (non necessarily continuous) then $\dim A=\dim B$.

\noindent $\bullet$ If $A_1,\dotsc, A_n$ are definable sets then
\[
\dim(A_1\cup\cdots\cup A_n)=\max\bigl\{\,\dim A_k;\;\;1\leq k\leq
n\,\bigr\}.
\]

\noindent $\bullet$ (\emph{Definable triviality of  semialgebraic
maps.}) We say that a semialgebraic map $\Phi: X\ra S$ is
\emph{definably trivial}   if there exists a    definable set
$F$, and a definable homeomorphism $\tau: X\ra F\times S$ such
that the diagram below is commutative
\[
\begin{diagram}
\node{X}\arrow{se,b}{\Phi}\arrow[2]{e,t}{\tau}\node[2]{S\times F}\arrow{sw,b}{\pi_S}\\
\node{}\node{S}\node{}
\end{diagram}.
\]
If $\Psi: X\ra Y$ is a continuous  definable map, and $p$ is a
positive integer, then there exists  a partition of $Y$ into
definable $C^p$-manifolds  $Y_1,\dotsc, Y_k$ such that     each
the restrictions
\[
\Psi: \Psi^{-1}(Y_k)\ra Y_k
\]
is definably trivial.

From the definable triviality  of  semialgebraic maps we deduce
the following consequence.

\begin{corollary} If $F: A\ra B$ is a continuous semialgebraic map, then the collection of fibers $(F^{-1}(b))_{b\in B}$ contains only finitely many homeomorphism types.\qed
\label{cor: homeo}
\end{corollary}

\begin{lemma} The discriminant set  $\Delta_\Gamma$ is a closed semialgebraic subset of  $S^2$ of dimension $\leq 1$.
\label{lemma: bertini}
\end{lemma}

\proof  We denote by $F(\Gamma)$ the set of points on the open
arcs where the curvature vanishes. The set $F(\Gamma)$ is       of
dimension  at most $1$, and thus is a finite  disjoint union of
(open arcs)  arcs and points. Moreover the collection
\[
\eL_\Gamma:=\bigl\{ L_p;\;\; p\in F(\Gamma)\,\bigr\}
\]
consists of finitely many affine lines. We see that
\[
\eC(\Gamma):=\Gamma\setminus\bigl(\, V\cup F(\Gamma)\,\bigr).
\]
The set $\Delta_\Gamma$ of degenerate vectors  is a disjoint union
\[
\Delta_\Gamma= \Delta_\Gamma^V\sqcup \Delta_\Gamma^F\sqcup
\Delta_\Gamma^*
\]

\begin{itemize}

\item $\Delta_\Gamma^V$ consists of unit vectors $\bu$
perpendicular to some  vector $\btau$ in some $\bsN_v$, $v\in V$.
It is a finite union of great circles on $S^2$

\item $\Delta_\Gamma^F$ consists of unit vectors $\bu$
perpendicular to one of the affine lines in $\eL_\Gamma$. It is
also    a finite union of great circles.

\item $\Delta_\Gamma^*$ consists   of unit vectors $\bu$ such that
there exist  a point $p\in \eC(\Gamma)$  with the property
$\bu\perp \bn(p),L_p$. It is $1$-dimensional and thus it is a
finite union of  semialgebraic arcs on $S^2$.

\end{itemize}

Observe that  for $p\in \eC(\Gamma)$ the subspace  spanned by
$\bn(p)$ and the tangent space $T_pe$  is the osculator plane of
$E$ at $p$. We denote this plane by $\eO_p$. Then
\[
\Delta^*_\Gamma=\bigl\{\bu\in S^2;\;\;\exists p\in
\eC(\Gamma):\;\;\bu\perp \eO_p\,\bigr\}.
\]
We need to prove that $\Delta^*_\Gamma$ is  has dimension $\leq
1$.  Consider the set
\[
\bsN^*(\Gamma)= \bigl\{ (p,u)\in \eC(\Gamma)\times S^2;\;\;
\bu\perp \eO_p\,\bigr\}.
\]
The set $\bsN^*(\Gamma)$ is semialgebraic,   and the natural map
$\bsN^*(\Gamma)\ra \eC(\Gamma)$ is two-to-one.  In particular
$\bsN^*(\Gamma)$ is one dimensional. The projection of
$\bsN^*(\Gamma)$ on $S^2$ is the set  $\Delta^*_\Gamma$ so that
$\dim \Delta^*_\Gamma\leq 1$. \qed

\noindent {\bf Proof of Proposition \ref{prop: morsesalg}.}
Suppose $\Gamma$ is a  compact, connected semialgebraic subset of
$\bR^3$ of dimension $1$. We define
\[
\eX_\Gamma:=\bigl\{ (\bu, p, q,\ve)\in S^2\times \Gamma\times
\Gamma
\times\bR;\;\;h_{\bu}(q)<h_{\bu}(p),\;\;|p-q|=\ve,\,\bigr\}.
\]
The natural projection $\pi:\eX_\Gamma\ra S^2\times \Gamma\bR$,
$(\bu,p,q,\ve)\mapsto (\bu,p,\ve)$ is a  continuous semialgebraic
map and its fiber  over $(\bu,p,\ve)$ is the  link
$L^+_\ve(p,\bu)$ defined in (\ref{eq: link}). From  Corollary
\ref{cor: homeo}  we deduce again that the are only finitely many
topological types amongst these spaces.

Fix a graph structure on $\Gamma$ and consider the definable
subset $\Cr\subset (S^2\setminus \Delta_\Gamma)\times \Gamma$
\[
\Cr:=\bigl\{ (\bu,p)\in (S^2\setminus \Delta_\Gamma)\times
\Gamma;\;\; p\in\Cr(\bu)\,\bigr\}.
\]
The natural projection $\Cr\ra S^2\setminus \Delta_\Gamma$ is
semialgebraic  and the fiber   over $\bu\in S^2\setminus
\Delta_\Gamma$ is the critical set $\Cr(\bu)$.  Corollary
\ref{cor: homeo}  now implies that there are finitely many
topological types  in the family $\Cr(\bu)$, $u\in\bS^2\setminus
\Delta_\Gamma$.  Since the sets $\Cr(\bu)$ are finite sets we
deduce
\[
\sup\bigl\{ |\Cr(\bu)|;\;\;u\in S^2\setminus
\Delta_\Gamma\,\bigr\} >\infty.
\]
Now observe that
\[
M_{\bu}(t)=\lim_{\ve\searrow 0} \sum_{p\in \Cr(\bu)} P_{\ve,
p,\bu}(t),
\]
where $P_{\ve,\bu}(t)$ is the Poincar\'{e} polynomial of the
topological pair $(\,{\rm cone}\,(\,L^-_\ve(t,p)\,),
L^-_\ve(t,p)\,)$.  The map
\[
\Cr\times \Gamma\times \bR\ni(\bu,p,\ve)\mapsto P_{\ve,
p,\bu}(t)\in \bZ[t]
\]
is definable. In particular, its range is finite, and its level
sets are   semialgebraic sets. This implies that $w(\bu)$ is
semialgebraic. \qed

\begin{remark}   Denote by $\eS_\Gamma\subset S^2$ the set of $\Gamma$-nondegenerate unit vectors $\bu$ such that the restriction of $h_{\bu}$ to $\Cr(\bu)$ is injective.  The set $\eS_\Gamma$  is open and semi-algebraic, and the functions  corresponding to $\bu\in\eS_\Gamma$ are the so called  stable stratified Morse functions.  The results of R. Pignoni \cite{Pig} imply that the  function
\[
\eS_\Gamma\ni \bu\mapsto M_{\bu}(t)\in \bZ
\]
is constant on the connected components of $\eS_\Gamma$. \qed
\end{remark}

\end{document}